\documentclass[twoside]{article}
\usepackage{mor}

\received{}                              %
\revised{}                               %
\pubyear{0x}                             %
\pubmonth{Xxxxxxx}                       %
\volume{xx}                              %
\issue{x}                                %
\pages{xxx--xxx}                         %
\firstpage{0xxx}                         %
\DOI{10.1287/moor.xxxx.xxxx}             %
\startpagenumber{1}                      %

\title{Integer Polynomial Optimization in Fixed Dimension}
\ShortTitle{Integer Polynomial Optimization in Fixed Dimension}
\ShortAuthors{J.~A.~De Loera, R.~Hemmecke, M.~K\"oppe, R.~Weismantel}
\NumberOfAuthors{4}
\FirstAuthor{Jes\'us A. De Loera}
\FirstAuthorAddress{University of California, Dept. of Mathematics,
  Davis CA 95616, USA}
\FirstAuthorEmail{deloera@math.ucdavis.edu}
\FirstAuthorURL{http://www.math.ucdavis.edu/~deloera}
\SecondAuthor{Raymond Hemmecke}
\SecondAuthorAddress{Otto-von-Guericke-Universit\"at Magdeburg, FMA/IMO,
  Universit\"atsplatz 2, 39106 Magdeburg, Germany}
\SecondAuthorEmail{hemmecke@imo.math.uni-magdeburg.de}
\SecondAuthorURL{http://www.math.uni-magdeburg.de/~hemmecke}
\ThirdAuthor{Matthias K\"oppe}
\ThirdAuthorAddress{Otto-von-Guericke-Universit\"at Magdeburg, FMA/IMO,
  Universit\"atsplatz 2, 39106 Magdeburg, Germany}
\ThirdAuthorEmail{mkoeppe@imo.math.uni-magdeburg.de}
\ThirdAuthorURL{http://www.math.uni-magdeburg.de/~mkoeppe}
\FourthAuthor{Robert Weismantel}
\FourthAuthorAddress{Otto-von-Guericke-Universit\"at Magdeburg, FMA/IMO,
  Universit\"atsplatz 2, 39106 Magdeburg, Germany}
\FourthAuthorEmail{weismant@imo.math.uni-magdeburg.de}
\FourthAuthorURL{http://www.math.uni-magdeburg.de/~weismant}
\FifthAuthor{}
\FifthAuthorAddress{}
\FifthAuthorEmail{}
\FifthAuthorURL{}


\MSCcodes{Primary: 90C10, 90C30; Secondary: 90C60, 90C57  } 

\ORMScodes{Primary:  Programming, Integer, Nonlinear; Secondary: Analysis of algorithms, Computational complexity} 

\usepackage{amssymb}
\usepackage{amsmath}
\usepackage{epsfig}
\usepackage{booktabs}
\usepackage{url}

\newcommand{\R}{\mathbb R}

\newcommand{\Z}{\mathbb Z}

\begin{document}
\maketitle

\begin{abstract}
  We classify, according to their computational complexity, integer
  optimization problems whose constraints and objective functions are
  polynomials with integer coefficients and the number of variables is
  fixed. For the optimization of an integer polynomial over 
  the lattice points of a convex polytope, we show an algorithm to compute 
  lower and upper bounds for the optimal value. 
  For polynomials that are non-negative over the polytope,
  these sequences of bounds lead to a fully polynomial-time approximation
  scheme for the optimization problem.
\end{abstract}

\section{Introduction}

Mixed integer non-linear programs combine the hardness of
combinatorial explosion with the non-convexity of non-linear
functions. For example, the well-known optimality conditions
developed for differentiable objective functions have no meaning when
the variables are discrete. Thus, it is perhaps not surprising that
already linear integer programming with general quadratic constraints
is undecidable \cite{jeroslow}. Nevertheless, when the number of
variables is fixed discrete optimization problems often become
tractable and efficient polynomial algorithms exist (e.g.\
\cite{bar,khachiyanporkolab,lenstra}). It is thus natural to ask {\em
  what is the complexity of integer non-linear optimization assuming
  that the number of variables is fixed?} We study the problem
\begin{equation} 
\label{model} \hbox{maximize} \  f(x_1,\dots,x_d) \
  \hbox{subject to} \ g_i(x_1,\dots,x_d) \geq 0, \ x \in  \Z^d. 
\end{equation}
Here $f,g_i$ are polynomials with integral coefficients. Note that
all throughout the paper we assume that the number of variables is
fixed. Here are our two contributions to the theory:

\noindent \emph{(1)} We give a classification of the computational
complexity of Problem (\ref{model}) according to special
cases. Section \ref{bounds} of this article presents the details, but
the reader can see the classification in Table \ref{gogo}. New results
are marked with letters, known results are marked with asterisks, arrows
indicate implications:

\begin{table}[htbp]
\small
\caption{Computational complexity of problem (\ref{model}) in fixed dimension.}
\label{gogo}
\begin{center}
\begin{tabular}{cc@{ }c@{ }c@{ }c@{ }c} \toprule
& \multicolumn{5}{c}{Type of objective function}\\
\cmidrule{2-6}
               &        & & convex      & & arbitrary  \\ 
Type of constraints & linear           & & polynomial  & & polynomial  \\ \midrule
Linear constraints, integer variables               & polytime ($*$) & $\Leftarrow$ &  polytime ($**$) &  & NP-hard (a)   \\ 
& $\Uparrow$ & & $\Uparrow$ & & $\Downarrow$ \\
Convex semialgebraic constraints, integer variables & polytime ($**$) & $\Leftarrow$ & polytime ($**$) &  & NP-hard (c)   \\
\\
Arbitrary polynomial constraints, integer variables & undecidable (b) & $\Rightarrow$ & undecidable (d) & $\Rightarrow$ & undecidable (e) \\
\bottomrule
\end{tabular}
\end{center}
\end{table}

\noindent {\em (2)} For problem $(a)$, that of optimizing an arbitrary 
integral polynomial over the lattice points of a convex rational
polytope with fixed number of variables, we present an algorithm to
compute a sequence of upper and lower bounds for its optimal value.
Our bounds can be used, for instance, in a branch-and-bound search for
the optimum.  We use Barvinok's algebraic encoding of the lattice points of
polytopes via rational functions \cite{BarviPom}. In Section \ref{algo}
we prove:

\begin{theorem} \label{main} Let the number of variables $d$ be fixed.  
  Let $f(x_1,\dots,x_d)$ be a polynomial of maximum total degree $D$ with
  integer coefficients, and let $P$ be a convex rational polytope
  defined by linear inequalities in $d$ variables.  We obtain an
  increasing sequence of lower bounds $\{L_k\}$ and a decreasing
  sequence of upper bounds $\{U_k\}$ to the optimal value
  \begin{equation}
    f^*=\hbox{maximize} \  f(x_1,x_2,\dots,x_d) \ \hbox{subject to} \ x
    \in P \cap \Z^d.\label{eq:maximize}
  \end{equation}
  The bounds $L_k$, $U_k$ can be computed in time polynomial in~$k$,
  the input size of $P$ and $f$, and the maximum total degree~$D$ and
  they satisfy the inequality $U_k-L_k \leq f^* \cdot (\sqrt[k]{|P
    \cap \Z^d|}-1).$
  
  More strongly, if $f$ is non-negative over the polytope (i.e.
  $f(x)\geq 0$ for all $x \in P$), there exists a fully
  polynomial-time approximation scheme (FPTAS) for the optimization
  problem~\eqref{eq:maximize}.

\end{theorem} 
\noindent We conclude with examples and a brief look at the mixed integer problem.
\section{Computational Complexity Bounds} \label{bounds}

All the results we present refer to the complexity model where the
number of operations is given in terms of the input size measured in
the standard binary encoding.  The results of H.~W. Lenstra Jr.
\cite{lenstra} imply the entry of Table \ref{gogo} marked with $(*)$,
i.e. solving linear integer programming problems with a fixed number
of variables can be done in time polynomial in the size of the input.
More recently, Khachiyan and Porkolab \cite{khachiyanporkolab} have
proved that in fixed dimension, the problem of minimizing a convex
polynomial objective function over the integers, subject to polynomial
constraints that define a convex body, can be solved in polynomial
time in the encoding length of the input.  Thus, they settled all
entries marked by $(**)$. By the natural containment exhibited by
these complexity classes, to show the validity of the remaining
entries of Table \ref{gogo} is enough to prove the following lemma:

\begin{lemma} 
\begin{enumerate}
\item The problem of minimizing a degree four polynomial over the lattice
  points of a convex polygon is NP-hard (entry $(a)$ in Table \ref{gogo}).

\item The problem of minimizing a linear form over polynomial
  constraints in at most 10 integer variables is not computable by a
  recursive function (entry $(b)$ in Table \ref{gogo}).
\end{enumerate}
\end{lemma}

\begin{proof} (1) We use the NP-complete problem AN1 on page 249 of
\cite{GarJohn79}. This problem states it is NP-complete to decide whether,
given three positive integers $a,b,c$, there exists a positive integer
$x<c$ such that $x^2$ is congruent with $a$ modulo $b$.  This problem is
clearly equivalent to asking whether the minimum of the quartic
polynomial function $(x^2-a-by)^2$ over the lattice points of 
the rectangle $\{(x,y) | \ 1 \leq x  \leq c-1,\ \frac{1-a}{b} \leq y
\leq \frac{(c-1)^2-a}{b} \}$ is zero or not. This settles part (1).

(2) In 1973 Jeroslow \cite{jeroslow} proved a similar result 
without fixing the number of variables. We follow his idea, but
resorting to a stronger lemma. More precisely our proof relies on a
1982 result \cite{jones} which states that there is no recursive
function that, given an integer polynomial $f$ with \emph{nine}
variables, can determine whether $f$ has a non-negative integer
zero, in the sense that it finds an explicit zero or returns
null otherwise.  Jones paper is a strengthening of the original
solution of Hilbert's tenth problem \cite{mati}. Now to each
polynomial $f$ in $\Z [x_1,x_2\dots,x_9]$ associate the
ten-dimensional minimization problem
\begin{equation} 
\hbox{minimize} \ y  \quad  \hbox
{subject to} \quad (1-y)f(x_1,x_2,\dots,x_9)=0, \ (y,x_1,\dots,x_9)
\in \Z_{\geq 0}^{10}. 
\end{equation}
The minimum attained by $y$ is either zero or one depending on whether
$f$ has an integer non-negative solution or not. Thus part (2) is
settled.
\end{proof}

\section{FPTAS for Optimizing Non-Negative Polynomials over
  Integer Points of Polytopes} \label{algo}

Consider now a polynomial function $f \in
\Z[x_1,x_2, \dots ,x_d]$ of maximum total degree $D$ and a
convex polytope $P=\{x | Ax \leq b\}$ where $A$ is an $m \times
d$ integral matrix and $b$ is an integral $m$-vector. The purpose of this
section is to present an algorithm to generate lower and upper bounds
$L_k,U_k$ to the integer global optimum value of
\begin{equation} 
\hbox{maximize} \ f(x_1,\dots,x_d) \  \hbox
{subject to} \  (x_1,\dots,x_d) \in P \cap \Z^d. 
\end{equation}
We should also remark that in our algorithm the polynomial objective
function $f$ can be arbitrary (e.g.\ non-convex). As we have seen, the
optimization problem is NP-hard already for two integer variables and
polynomials of degree four. Nevertheless we will see that, in fixed
dimension and when $f(x) \geq 0$ for all $x \in P$, the algorithm gives a
fully polynomial time approximation scheme or FPTAS. This means that,
in polynomial time on the input and $(1/\epsilon)$, one can compute a
$(1-\epsilon)$-approximation to the maximum. The algorithm we present
is based on A. Barvinok's theory for encoding all the lattice points
of a polyhedron in terms of short rational functions. See
\cite{bar,BarviPom} for all details.  Lattice points are thought of as
exponent vectors of monomials.  For example, $z_1^2z_2^{-11}$ encodes
the lattice point $(2,-11)$.  The set of lattice points is represented
by a Laurent polynomial: $g_P(z)=\sum_{\alpha \in P\cap\Z^d}
z^\alpha.$ From Barvinok's theory this exponentially-large sum of
monomials $g_P(z)$ can instead be written as a polynomial-size sum of
rational functions (assuming the dimension $d$ is fixed) of the form:
\begin{equation} \label{eq:aa}
g_P(z) = \sum_{i\in I} {E_i \frac{z^{u_i}} {\prod\limits_{j=1}^d
(1-z^{v_{ij}})}},
\end{equation}
where $I$ is a polynomial-size indexing set, and where
$E_i\in\{1,-1\}$ and $u_i, v_{ij} \in\Z^d$ for all $i$ and $j$. For
details see \cite{BarviPom,latte1,latte2}.

We need a way to encode via rational functions the values of the
polynomial $f$ over all the lattice points in a polytope.  The key
idea, first introduced in Lemma 9 of \cite{latte2} and generalized in
\cite{hugginsthesis}, is that differential operators associated to $f$
can be used to compute a rational function representation of
$\sum_{a\in P \cap \Z^d} f(a) z^a.$ The following Lemma recently
appeared in \cite{newbar}:

\begin{lemma} \label{operators} Let $g_P(z)$ be the Barvinok 
  representation of the generating function of the lattice points of
  $P$. Let $f$ be a polynomial in $\Z[x_1,\dots,x_d]$ of maximum total
  degree $D$. We can compute, in time polynomial on $D$ and the size of the
  input data, a Barvinok rational function representation $g_{P,f}(z)$
  for the generating function $\sum_{a\in P \cap \Z^d} f(a) z^a.$
\end{lemma}

\begin{proof} We give here the author's original proof the lemma for $D$
  fixed. The first proof without this assumption was recently given by
  A. Barvinok in \cite{newbar}. 
  
  We begin assuming $f(z)=z_r$, the general case will follow from it:
  Consider the action of the differential operator
  $z_r\frac{\partial}{\partial z_r}$ in the generating function
  $g_P(z)$ and on its Barvinok representation.  On one hand, for the
  generating function
\[
z_r\frac{\partial}{\partial z_r}\cdot g_P(z)=\sum\limits_{\alpha\in
  P \cap \Z^d} z_r \frac{\partial}{\partial z_r} z^\alpha=
\sum\limits_{\alpha\in P \cap \Z^d}\alpha_r z^\alpha.
\]
On the other hand, by linearity of the operator, we have that in terms
of rational functions
\[
z_r\frac{\partial}{\partial z_r}\cdot g_P(z)=\sum\limits_{i \in I} E_i
z_r\frac{\partial}{\partial z_r}\cdot \left( \frac{z^{u_i}} {\prod
\limits_{j=1}^d (1-z^{v_{ij}})}\right).
\]
Thus it is enough to prove that the summands of the expression above
can be written in terms of rational functions computable in polynomial
time. The standard quotient rule for derivatives says that
\[
\frac{\partial}{\partial z_r} \left(\frac{z^{u_i}} {\prod
\limits_{j=1}^d (1-z^{v_{ij}})}\right)=\frac{(\frac{\partial z^{u_i} }{\partial
z_r}) \prod^d_{j=1}(1-z^{v_{ij}})
-z^{u_i}(\frac{\partial}{\partial z_r}
\prod_{j=1}^d(1-z^{v_{ij}}))} {\prod_{j=1}^d(1-z^{v_{ij}})^2}.
\]
We can expand the numerator as a sum of no more than $2^d$ monomials.
This is a constant number because $d$, the number of variables, is
assumed to be a constant. This argument completes the proof of
our lemma when $f(z)=z_r$. 

For the case when $f(z)$ is a general monomial, i.e. $f(z)=c\cdot
z_1^{\beta_1}\cdot\ldots\cdot z_d^{\beta_d}$, then we can compute
again a rational function representation of $g_{P,f}(z)$ by repeated
application of basic differential operators:
\[
c\left(z_1\frac{\partial}{\partial z_1}\right)^{\beta_1}
\cdot\ldots\cdot \left(z_d \frac{\partial}{\partial
z_d}\right)^{\beta_d}\cdot g_P(z)=\sum\limits_{\alpha\in P \cap \Z^d} c\cdot
\alpha^{\beta} z^\alpha.
\]
Thus we require no more than $O(D^d)$ repetitions of the single-variable
case. 

Finally, if we deal with a polynomial $f$ of many monomial terms, we
compute and add up all such expressions that we get for each term of $f(x)$ and
obtain our desired short rational function representation for the
generating function for $ \sum_{\alpha\in P \cap \Z^d} f(\alpha)
z^\alpha.$ Note that only polynomially many steps are needed because
$d$ is fixed and the largest number of possible monomials in $f$ 
of degree $s$ is $\left(\begin{smallmatrix}d+s-1\\ d-1\\
  \end{smallmatrix}\right)$, thus for fixed $d$ we will do no more than
$O(D^{d})$ repetitions of the monomial case.
\end{proof}

Now we are ready to present our algorithm to obtain bounds $U_k,L_k$
that reach the optimum.  Step 1 of preprocessing is necessary because
we rely on the elementary fact that, for a collection
$S=\{s_1,\dots,s_r\}$ of non-negative real numbers, $\hbox{\rm
  maximum} \{s_i | s_i \in S \}$ equals $\lim_{k \rightarrow \infty}
\sqrt[k]{\sum_{j=1}^r s_j^k}.$

\noindent {{\bf Algorithm}} \label{Algorithm}

\noindent {\em Input:} A rational convex polytope $P \subset \R^d$, 
a polynomial objective $f \in \Z[x_1,\dots,x_d]$ of maximum total degree $D$.

\noindent {\em Output:} An increasing sequence of lower bounds $L_k$,
and a decreasing sequence of upper bounds $U_k$ reaching the maximal
function value $f^*$ of $f$ over all lattice points of $P$.

\noindent {\bf Step 1.} If $f$ is known to be 
non-negative in all points of $P$, then go directly to Step 2. 
Else, solving $2d$ linear programs over $P$, we
find lower and upper integer bounds for each of the variables
$x_1,\ldots,x_d$. Let $M$ be the maximum of the absolute values of
these $2d$ numbers. Thus $|x_i|\leq M$ for all $i$. Let $C$ be the
maximum of the absolute values of all coefficients, and $r$ be the
number of monomials of $f(x)$.  Then
\[
L:=-rCM^D\leq f(x)\leq rCM^D=:U,
\]
as we can bound the absolute value of each monomial of $f(x)$ by
$CM^D$. Replace $f$ by $\bar{f}(x)=f(x)-L\leq U-L$, a non-negative
polynomial over $P$. Go to Steps 2, 3, etc.~and return the optimal
value of $\bar{f}.$ Trivially, if we find the optimal value of
$\bar{f}$ over $P$ we can extract the optimal value for $f$.

\noindent {\bf Step 2.} Via Barvinok's algorithm (see
\cite{bar,BarviPom,newbar}), compute a short rational function
expression for the generating function $g_P(z)=\sum_{\alpha\in
  P\cap\Z^d} z^\alpha$.  From $g_P(z)$ compute the number $|P \cap
\Z^d|=g_P(1)$ of lattice points in $P$ in polynomial time.

\noindent {\bf Step 3.} From the rational function representation $ g_P(z)$ of the
generating function $\sum\limits_{\alpha\in P \cap \Z^d}
z^\alpha$ compute the rational function representation of
$g_{P,f^k}(z)$ of $\sum_{\alpha\in P \cap \Z^d} f^k(\alpha) z^\alpha$ in
polynomial time by application of Lemma \ref{operators}. We define

\[
L_k:=\sqrt[k]{g_{P,f^k}(1)/g_{P,f^0}(1)}\;\;\;\text{and}\;\;\;
U_k:=\sqrt[k]{g_{P,f^k}(1)}.
\]

When $\lfloor U_{k}\rfloor-\lceil L_{k}\rceil<1$  stop and return
$\lceil L_{k} \rceil=\lfloor U_{k} \rfloor$ as the optimal value.

\noindent {\bf End of Algorithm}.

\begin{lemma} The algorithm is correct. \end{lemma}

\begin{proof} Using the fact that the arithmetic mean of a finite set
  of nonnegative values is at most as big as the maximum value, which
  in turn is at most as big as the sum of all values, we obtain the
  sequences of lower and upper bounds, $L_{k}$ and $U_{k}$, for
  the maximum:
\[
L_{k}=
\sqrt[k]{\frac{\sum\limits_{\alpha\in P \cap \Z^d}
f(\alpha)^{k}}{|P\cap\Z^d|}} \leq
\max\{f(\alpha):\alpha\in P\cap\Z^d\}\leq
\sqrt[k]{\sum\limits_{\alpha\in P \cap \Z^d}
f(\alpha)^{k}}=U_{k}.
\]
Note that as $s\rightarrow\infty$, $L_{k}$ and $U_{k}$ approach
this maximum value monotonously (from below and above, respectively).
Trivially, if the difference between (rounded) upper and lower bounds
becomes strictly less than $1$, we have determined the value
$\max\{f(x):x\in P\cap\Z^d\}=\lceil L_{k}\rceil$. Thus the algorithm
terminates with the correct answer.
\end{proof}

Theorem \ref{main} will follow from the next lemma:

\begin{lemma}\label{lemma:bounds}
  Let $f$ be a polynomial with integer coefficients and maximum total
  degree~$D$. When the dimension $d$ is fixed,

\begin{enumerate}  
\item the bounds $L_k$, $U_k$ can be computed in time polynomial in
  $k$, the input size of $P$ and $f$, and the total degree~$D$. The bounds
  satisfy the following inequality:
$$ 
U_k-L_k \leq f^* \cdot \left(\sqrt[k]{|P \cap \Z^d|}-1 \right).
$$
\item In addition, when $f$ is non-negative over $P$ (i.e. $f(x)\geq 0$
  for all $x \in P$), for $k=(1+1/\epsilon)\log({|P \cap \Z^d|})$,
  $L_k$ is a $(1-\epsilon)$-approximation to the optimal value $f^*$ and it
  can be computed in time polynomial in the input size, the total
  degree~$D$, and $1/\epsilon$. Similarly, $U_k$ gives a
  $(1+\epsilon)$-approximation to $f^*$. Moreover, with the same
  complexity, one can also find a feasible lattice point that
  approximates an optimal solution with similar quality.
\end{enumerate}
\end{lemma}
\begin{proof}
\noindent  Part (i). From Lemma \ref{operators} on fixed dimension $d$, we can
compute $g_{P,f}=\sum_{\alpha\in P \cap \Z^d} f(\alpha) z^\alpha$ as a
rational function in time polynomial in $D$, the total degree of $f$,
and the input size of $P$. Thus, because $f^k$ has total degree
of~$Dk$ and the encoding length for the coefficients of $f^k$ is
bounded by $k \log (kC)$ (with $C$ the largest coefficient in $f$), we can also compute $g_{P,f^k} = \sum_{\alpha\in P \cap \Z^d}
f^k(\alpha) z^\alpha$ in time polynomial in $k$, the total degree~$D$,
and the input size of~$P$.  Note that using residue techniques
\cite{newbar}, we can evaluate $g_{P,f^k}(1)$ in polynomial time.
Finally observe

\begin{align*}
U_k-L_k&= \sqrt[k]{\sum_{\alpha \in P \cap \Z^d}
  f^k(\alpha)}-\sqrt[k]{\frac{\sum_{\alpha \in P \cap \Z^d}
  f^k(\alpha)}{|P \cap \Z^d|}} = \sqrt[k]{\frac{\sum_{\alpha \in P
  \cap \Z^d} f^k(\alpha)}{|P \cap \Z^d|}} \left( \sqrt[k]{|P \cap
  \Z^d|}-1 \right) \\
& =L_k \left( \sqrt[k]{|P \cap \Z^d|}-1 \right) \leq f^*
  \left(\sqrt[k]{|P \cap \Z^d|}-1 \right). 
\end{align*}

\noindent Part (ii). Note that if $\left(\sqrt[k]{|P \cap \Z^d|}-1 \right) \leq
\epsilon$ then $L_k$ is indeed a $(1-\epsilon)$-approximation because 
$$
f^* \leq U_k =L_k+(U_k-L_k) \leq L_k + f^*\left(\sqrt[k]{|P \cap
    \Z^d|}-1 \right) \leq L_k +f^* \epsilon.$$

Observe that $\phi(\epsilon) := (1+1/\epsilon)/(1/\log(1+\epsilon))$ is an
increasing function for $\epsilon <1$ and $\lim_{\epsilon \rightarrow
  0}\phi(\epsilon)=1$, thus $\phi(\epsilon)\geq 1$ for $0<\epsilon\leq1$. 
Hence, for all 
$k \geq \log({|P \cap \Z^d|})+\log({|P \cap
  \Z^d|})/\epsilon \geq \log({|P \cap \Z^d|})/\log(1+\epsilon),$ 
we have indeed $\left(\sqrt[k]{|P \cap \Z^d|}-1 \right) \leq \epsilon.$
Finally, from Lemma \ref{operators}, the calculation of $L_k$ for
$k=\log({|P \cap \Z^d|})+\log({|P \cap \Z^d|})/\epsilon$ would require a
number of steps polynomial in the input size and $1/\epsilon$. A
very similar argument can be written for $U_k$ but we omit it here.

%

To complete the proof of part (ii) it remains to show that not only
we approximate the optimal value $f^*$ but we can also efficiently
find a lattice point $\alpha$ with $f(\alpha)$ giving that quality
approximation of $f^*$. Let $k = (1+1/\epsilon)\log({|P \cap
  \Z^d|})$, thus, by the above discussion, $L_k$ is an
$(1-\epsilon)$-approximation to~$f^*$.  Let $Q_0 := [-M,M]^d$ denote
the box computed in Step~1 of the algorithm such that $P\subseteq
Q_0$.  By bisecting $Q_0$, we obtain two boxes $Q_1'$ and $Q_1''$.  By
applying the algorithm separately to the polyhedra $P\cap Q_1'$ and
$P\cap Q_1''$, we compute lower bounds $L_k'$ and $L_k''$ for the
optimization problems restricted to $Q_1'$ and $Q_1''$, respectively.
Because $L_k^k$ is the arithmetic mean of $f^k(\alpha)$ for $\alpha\in
P\cap\Z^d$, clearly
\begin{displaymath}
  \min\{L_k', L_k''\} \leq L_k \leq \max\{L_k', L_k''\}.
\end{displaymath}
Without loss of generality, let $L_k'\geq L_k''$.  We now apply the
bisection procedure iteratively on $Q_k'$.  After $d\log M$ bisection
steps, we obtain a box $Q_k'$ that contains a single lattice point
$\alpha\in P\cap Q_k'\cap Z^d$, which has an objective value
$f(\alpha) = L_k' \geq L_k \geq (1-\epsilon) f^*$.
\end{proof}

We remark that if we need to apply the construction of Step~1 of the
algorithm because $f$ takes negative values on~$P$, then we can only
obtain an $(1-\epsilon)$-approximation (and
$(1+\epsilon)$-approximation, respectively) for the modified function
$\bar f$ in polynomial time, but not the original function $f$. We
also emphasize that, although our algorithm requires the computation
of $\sum_{\alpha \in P} f^q(\alpha)$ for different powers of $f$,
these numbers are obtained without explicitly listing all lattice
points (a hard task), nor we assume any knowledge of the individual
values $f(\alpha)$. We can access the power means $\sum_{\alpha \in P}
f^q(\alpha)$ indirectly via rational functions. Here are two small
examples:

\paragraph{Example 1, monomial optimization over a quadrilateral:} 
The problem we consider is that of maximizing the value of the  
monomial $x^3y$ over the lattice points of the quadrilateral
$$
\{(x,y) | 3991 \leq 3996\,x-4\,y \leq 3993, \ 1/2\leq x \leq 5/2\}.
$$
It contains only 2 lattice points. The sum of rational functions
encoding the lattice points is
\[
\frac{{x}^{2}{y}^{1000}}{\left( 1-(xy^{999})^{-1} \right)
 \left( 1-{y}^{-1} \right)}+\frac{xy}{\left( 1-x{y}^{999} \right)
 \left( 1-{y}^{-1} \right)}+
{\frac {xy}{ \left( 1-x{y}^{999}
 \right)  \left( 1-y \right) }}+\frac{{x}^{2}{y}^{1000}}{\left( 1-(xy^{999})^{-1} \right)\left( 1-y \right)}.
\]
In the first iteration $L_1=4000.50$ while $U_1=8001$. After thirty iterations,
we see $L_{30}=7817.279750$ while $U_{30}=8000$, the true optimal value.

\paragraph{Example 2, nvs04 from {\small MINLPLIB}:}
A somewhat more complicated example, from a well-known
library of test examples (see \url{http://www.gamsworld.org/minlp/}),
is the problem given by
\begin{equation}
  \begin{aligned}
    \min\quad & 100 \left(\frac12 + i_2 - \left(\frac35 + i_1\right)^2 \right)^2 + \left(\frac25 -
    i_1\right)^2\\
    \text{s.\,t.}\quad & i_1, i_2 \in [0,200] \cap \Z.
  \end{aligned}
\end{equation}
Its optimal solution as given in {\small MINLPLIB} is $i_1=1$, $i_2=2$ with an
objective value of~$0.72$. Clearly, to apply our algorithm from page
\pageref{Algorithm} literally, the objective function needs to be
multiplied by a factor of $100$ to obtain an integer valued polynomial.

Using the bounds on $i_1$ and $i_2$ we obtain an upper bound of $165
\cdot 10^9$ for the objective function, which allows us to convert the
problem into an equivalent maximization problem, where all feasible
points have a non-negative objective value. The new optimal objective
value is $164999999999.28.$ Expanding the new
objective function and translating it into a differential operator
yields
  \begin{align*}
    &\frac{4124999999947}{25} \mathrm{Id} - 28 z_2 \frac{\partial}{\partial
      z_2}+\frac{172}{5} z_1\frac{\partial}{\partial z_1}-117 \left(z_1
      \frac{\partial}{\partial z_1} \right) ^{(2)} - 100
    \left(z_2 \frac{\partial}{\partial z_2}\right)^{(2)} \\
    &+ 240 \left(z_2
      \frac{\partial}{\partial z_2}\right) \left(z_1 \frac{\partial}{\partial
        z_1}\right) 
    + 200  \left(z_2
      \frac{\partial}{\partial z_2}\right) \left(z_1 \frac{\partial}{\partial
        z_1}\right) ^{(2)} - 240 \left(z_1 \frac{\partial}{\partial
        z_1}\right)^{(3)}  
    - 100 \left(z_1 \frac{\partial}{\partial
        z_1}\right)^{(4)}.
  \end{align*}
The short generating function can be written as 
$ g(z_1,z_2) = \left( \frac1{1-z_1} - \frac{z_1^{201}}{1-z_1} \right)
  \left( \frac1{1-z_2} - \frac{z_2^{201}}{1-z_2} \right).$

In this example, the number of lattice points is $|P\cap \Z^2|=40401.$
The first bounds are $L_1=139463892042.292155534$,
$U_1=28032242300500.723262442$.  After 30 iterations the bounds become
$L_{30}=164999998845.993553019$ and $U_{30}=165000000475.892451381.$

\section{An extension to the mixed integer case}
Now, we wish to discuss extensions of Theorem \ref{main} to the
mixed integer scenario. If some of the variables are allowed to be
continuous then we can describe the task as
$$
\hbox{maximize} \ f(x,y) \ \hbox{subject to} \ \{ (x,y) | Ax+By \leq
b, \ \hbox{with} \ x_i \in \Z \ \hbox{and} \ y_i \in \R\}.
$$
Consider the sequence of integer polynomial optimization problems
$\hbox{optimize} \ f(x,\frac{y}{n}) \ \hbox{subject to}\ (x,y) \in
\Gamma_n=\{ (x,y) | Ax+B(\frac{y}{n}) \leq b, \ \hbox{with} \ x_i \in
\Z \ \hbox{and} \ y_i \in \Z \},$ where each of the subproblems is
equivalent to optimizing a polynomial over a ``semi-dilated'' polytope
(in some coordinate directions but not others). As $n$ goes to
infinity, the sequence of optimal solution values can
have several limit points.  Nevertheless, it is still possible 
to construct a subsequence 
of problems whose optimal values approximate the mixed integer optimum 
to arbitrary precision:

\begin{corollary} \label{mip} With the hypotheses of Theorem
  \ref{main}, 
  we can construct a sequence of integer polynomial programming problems, over
  finer and finer grids, whose optimal values converge
  to the optimal value of the mixed integer program
  \begin{displaymath}
    f^*=\hbox{maximize} \  f(x_1,x_2,\dots,x_d) \ \hbox{subject to} \ x \in P, \ \hbox{and} \ x_i \in \Z \
    \hbox{for} \ i \in I \subseteq  \{1,2,\dots,d\}.
  \end{displaymath}
  By applying the algorithm of Theorem~\ref{main} to the subproblems, we can approximate the
  optimum to arbitrary precision.
\end{corollary}

When all variables are continuous, the original
polytope $P$ is actually dilated uniformly in all directions by a
parameter $n$, or equivalently, the integer grid is refined. Define
$$ 
r(q,n)=\sum_{\alpha \in (nP) \cap \Z^d} f^q(\alpha) z^\alpha,
$$
for each power $q$ and dilation factor $n$. Note that for fixed
$q$, then we can easily see that the sum $\sum_{\alpha \in (nP) \cap
  \Z^d} \frac{1}{n^d} f^q(\alpha) z^\alpha$ is essentially an
approximation to the Riemann integral of $f^q$; thus
$$ 
\lim_{n \rightarrow \infty} r(q,n)/n^d = \int_P f^q(\bar{x})
d\bar{x}.
$$ 
As $n,q$ grow, the values $\sqrt[q]{r(q,n)}$ approximate the sequence
$\sqrt[q]{\int_P f^q(\bar{x}) d\bar{x}}$ which converges to
$(\hbox{max}_{x \in P} f(x))\cdot \hbox{volume} (P)$.  This is related
to recent work (see \cite{deklerketal,lasserre,parrilosturmfels} and
references therein) where the global optimum of a polynomial over a
compact domain is investigated as the result of a grid refinement and
properties of sums of squares.

\acknowledgments{We are truly grateful to Prof. Alexander Barvinok who
  communicated to us that Lemma \ref{operators} was true for variable
  $D$ and thus we had indeed obtained an FPTAS from the construction
  of the upper and lower bounds. We thank the anonymous referees for
  the many suggestions that improved the presentation. The first
  author gratefully acknowledges support from NSF grant DMS-0309694, a
  2003 UC-Davis Chancellor's fellow award, the Alexander von Humboldt
  foundation, and IMO-Magdeburg. The remaining authors were supported
  by the European TMR network ADONET 504438.}

\end{document}